\documentclass[12pt]{article}
\usepackage{amssymb, amsmath, amsthm, amscd}
\usepackage[dvips]{graphics}
\usepackage[utf8]{inputenc}
\usepackage[all,cmtip]{xy}
\usepackage{bbm}
\usepackage{enumitem}
\usepackage{setspace}
\usepackage[colorlinks=true]{hyperref}
\hypersetup{colorlinks   = true}
\hypersetup{linkcolor=blue}
\begin{document}

\newtheorem{The}{Theorem}[section]
\newtheorem{Lem}[The]{Lemma}
\newtheorem{Prop}[The]{Proposition}
\newtheorem{Cor}[The]{Corollary}
\newtheorem{Rem}[The]{Remark}
\newtheorem{Obs}[The]{Observation}
\newtheorem{SConj}[The]{Standard Conjecture}
\newtheorem{Titre}[The]{\!\!\!\! }
\newtheorem{Conj}[The]{Conjecture}
\newtheorem{Question}[The]{Question}
\newtheorem{Prob}[The]{Problem}
\newtheorem{Def}[The]{Definition}
\newtheorem{Not}[The]{Notation}
\newtheorem{Claim}[The]{Claim}
\newtheorem{Conc}[The]{Conclusion}
\newtheorem{Ex}[The]{Example}
\newtheorem{Fact}[The]{Fact}
\newtheorem{Formula}[The]{Formula}
\newtheorem{Formulae}[The]{Formulae}
\newtheorem{Prop-Def}[The]{Proposition and Definition}
\newcommand{\C}{\mathbb{C}}
\newcommand{\R}{\mathbb{R}}
\newcommand{\N}{\mathbb{N}}
\newcommand{\Z}{\mathbb{Z}}
\newcommand{\Q}{\mathbb{Q}}
\newcommand{\Proj}{\mathbb{P}}
\newcommand{\Rc}{\mathcal{R}}
\newcommand{\Oc}{\mathcal{O}}
\newcommand\george[1]{{\textcolor{green}{#1}}}
\newcommand\dan[1]{{\textcolor{blue}{#1}}}

\begin{center}

{\Large\bf Non-K\"ahler Hyperbolicity:\\[0.5ex] Surfaces and Codimension One}

\end{center}

\begin{center}

  {\large  Samir Marouani}\\[1ex]
  Universit\'e de Kairouan,\\
  Institut Sup\'erieur des \'Etudes Appliqu\'ees en Humanit\'es de Sbe\"itla, Tunisia\\[0.5ex]
  \texttt{marouanisamir20@gmail.com}

\end{center}

\begin{abstract}
We study hyperbolicity notions for compact complex manifolds defined by special Hermitian metrics and by holomorphic maps from $\mathbb C^{n-1}$ (divisorial hyperbolicity). For surfaces we show that divisorial and Brody hyperbolicity coincide and that Hermitian-symplectic and K\"ahler hyperbolicity coincide when $h^{2,0}=0$. In any dimension we prove obstructions to Gauduchon hyperbolicity, show that projective spaces are not divisorially hyperbolic, and give divisorially hyperbolic manifolds that are not balanced hyperbolic.
\end{abstract}

\vspace{3ex}
\section{Introduction}

\subsection*{Background}
A compact complex manifold is {\it Brody hyperbolic} if it contains no non-constant entire curve $\C\to X$ (equivalently, by Brody's theorem, Kobayashi hyperbolic), and {\it K\"ahler hyperbolic} in the sense of Gromov [Gro91] if it carries a K\"ahler metric whose lift to the universal cover has a bounded primitive. Several extensions to non-K\"ahler manifolds have been introduced:
\begin{itemize}
\item {\it balanced hyperbolicity} and {\it divisorial hyperbolicity} by Marouani--Popovici [MP22], where curves are replaced by divisors and $(1,1)$-forms by $(n-1,n-1)$-forms. We prove that balanced hyperbolic manifolds are divisorially hyperbolic, give examples (degenerate balanced manifolds, products), and introduce divisorially K\"ahler and divisorially nef classes together with several questions;
\item {\it SKT hyperbolicity} and {\it Gauduchon hyperbolicity} by Marouani [Mar24]: SKT hyperbolic manifolds are Kobayashi hyperbolic and Gauduchon hyperbolic manifolds are divisorially hyperbolic;
\item {\it Hermitian-symplectic (H-S) hyperbolicity} by Khelifati [Khe26], who proves that H-S hyperbolic manifolds are Kobayashi hyperbolic.
\end{itemize}

\subsection*{Implications between the notions}
In any dimension $n\geq2$ (Lemma~\ref{lem:impl}, [MP22], [Mar24]):
\begin{itemize}
\item K\"ahler hyperbolic $\Rightarrow$ balanced hyperbolic (equivalent if $n=2$) $\Rightarrow$ Gauduchon hyperbolic (equivalent to SKT hyperbolic if $n=2$) $\Rightarrow$ divisorially hyperbolic (equivalent to Brody hyperbolic if $n=2$, Proposition~\ref{prop:brody});
\item K\"ahler hyperbolic $\Rightarrow$ H-S hyperbolic $\Rightarrow$ SKT hyperbolic $\Rightarrow$ Kobayashi hyperbolic (equivalent to Brody hyperbolic) $\Rightarrow$ divisorially hyperbolic.
\end{itemize}

\subsection*{What is added and treated in this paper}
\begin{enumerate}[label=(\Alph*)]
\item {\bf Surfaces: divisorial and Brody hyperbolicity} (Section~\ref{sec:surfaces}). For surfaces, divisorial hyperbolicity coincides with Brody hyperbolicity (Proposition~\ref{prop:brody}).
\item {\bf Surfaces: divisorially K\"ahler classes and H-S hyperbolicity} (Section~\ref{sec:classes}). We describe explicitly the divisorially K\"ahler and nef cones of [MP22] on a K\"ahler surface and prove that their closures coincide (Proposition~\ref{prop:DK}, Corollary~\ref{cor:Q39}), which answers in dimension two a question of [MP22]. H-S hyperbolicity means that some $\widetilde d$-bounded class is divisorially K\"ahler (Theorem~\ref{thm:HSDK}); this gives the equivalence with K\"ahler hyperbolicity when $h^{2,0}=0$ (Corollary~\ref{cor:K}). Divisorially hyperbolic K\"ahler surfaces have a divisorially nef canonical class (Proposition~\ref{prop:Q15}), a second question of [MP22] treated in dimension two.
\item {\bf Obstructions in any dimension} (Section~\ref{sec:codim}). A Gauduchon hyperbolic manifold contains no compact hypersurface lifting to the universal cover (Proposition~\ref{prop:hyper}); in the K\"ahler case it has a non-amenable fundamental group (Lemma~\ref{lem:iso}, Corollary~\ref{cor:amen}); non-degenerate maps $\C^{n-1}\to X$ satisfy an isoperimetric inequality which forces exponential growth (Propositions~\ref{prop:iso} and~\ref{prop:exp}).
\item {\bf Examples} (Section~\ref{sec:ex}). No projective space $\mathbb P^n$ is divisorially hyperbolic (Proposition~\ref{prop:Pn}), a point left open in [MP22, \S2.3 (V)]. We give divisorially hyperbolic manifolds that are neither balanced nor Gauduchon hyperbolic in every dimension $n\geq2$, starting with surfaces (Table~\ref{tab:ex}).
\end{enumerate}
The few open questions we consider important are stated where they arise.

\section{Preliminaries and definitions}\label{sec:prelim}

Throughout, $X$ is a compact complex manifold of dimension $n\geq2$, $\omega$ is a Hermitian metric on $X$ (a $C^\infty$ positive definite $(1,1)$-form), $\pi:\widetilde X\to X$ is the universal cover and $\widetilde\omega=\pi^\star\omega$.

\subsection{Bounded forms}

A form $\varphi$ on $\widetilde X$ is {\bf bounded} if $\sup_{\widetilde X}|\varphi|_{\widetilde\omega}<\infty$. Since $X$ is compact, the lift $\pi^\star\varphi$ of every smooth form on $X$ is bounded, and boundedness does not depend on the choice of $\omega$ (two Hermitian metrics on $X$ are comparable).

\begin{Def}\label{def:bounded}
Let $\varphi$ be a $C^\infty$ $k$-form on $X$.
\begin{enumerate}[label=(\alph*)]
\item $\varphi$ is {\bf $\widetilde d$-bounded} if $\pi^\star\varphi=d\beta$ for a bounded $(k-1)$-form $\beta$ on $\widetilde X$ [Gro91].
\item $\varphi$ is {\bf $\widetilde{(\partial+\bar\partial)}$-bounded} if $\pi^\star\varphi=\partial\alpha+\bar\partial\beta$ for bounded $(k-1)$-forms $\alpha,\beta$ on $\widetilde X$ [Mar24].
\end{enumerate}
\end{Def}

\subsection{The hyperbolicity notions}

Recall the special Hermitian metrics: $\omega$ is {\it K\"ahler} if $d\omega=0$; {\it balanced} if $d\omega^{n-1}=0$; {\it SKT} (strong K\"ahler with torsion) if $\partial\bar\partial\omega=0$; {\it Gauduchon} if $\partial\bar\partial\omega^{n-1}=0$; {\it Hermitian-symplectic} (H-S) if there is a $(0,2)$-form $\rho$ with $d(\rho+\omega+\bar\rho)=0$. For $n=2$ one has $\omega^{n-1}=\omega$, so balanced $=$ K\"ahler and Gauduchon $=$ SKT; hence for surfaces balanced (resp.\ Gauduchon) hyperbolicity is K\"ahler (resp.\ SKT) hyperbolicity.

\begin{Def}\label{def:hyp}
The manifold $X$ is
\begin{enumerate}[label=(\alph*)]
\item {\bf K\"ahler hyperbolic} [Gro91] if it carries a K\"ahler metric $\omega$ that is $\widetilde d$-bounded;
\item {\bf balanced hyperbolic} [MP22] if it carries a balanced metric $\omega$ such that $\omega^{n-1}$ is $\widetilde d$-bounded;
\item {\bf SKT hyperbolic} [Mar24] if it carries an SKT metric $\omega$ that is $\widetilde{(\partial+\bar\partial)}$-bounded;
\item {\bf Gauduchon hyperbolic} [Mar24] if it carries a Hermitian metric $\omega$ such that $\omega^{n-1}$ is $\widetilde{(\partial+\bar\partial)}$-bounded;
\item {\bf H-S hyperbolic} [Khe26] if it carries a Hermitian metric $\omega$ and a $(0,2)$-form $\rho$ such that $\Omega=\rho+\omega+\bar\rho$ is $d$-closed and $\widetilde d$-bounded.
\end{enumerate}
\end{Def}

\begin{Lem}\label{lem:auto}
If $\pi^\star\omega=\partial\alpha+\bar\partial\beta$ for some smooth forms $\alpha,\beta$, then $\omega$ is SKT. If $\pi^\star\omega^{n-1}=\partial\Gamma+\bar\partial\Gamma'$ for some smooth forms $\Gamma,\Gamma'$, then $\omega$ is Gauduchon. In particular the SKT (resp.\ Gauduchon) condition in Definition~\ref{def:hyp}(c) (resp.\ (d)) is automatic.
\end{Lem}

\begin{proof}
Since $\bar\partial\partial=-\partial\bar\partial$, $\partial^2=0$ and $\bar\partial^2=0$,
\begin{eqnarray*}
\partial\bar\partial(\partial\alpha+\bar\partial\beta) & = & -\,\partial\partial\bar\partial\alpha\;+\;\partial\bar\partial\bar\partial\beta\;=\;0 .
\end{eqnarray*}
Hence $\partial\bar\partial\,\pi^\star\omega=0$. As $\pi$ is a local biholomorphism, $\partial\bar\partial\,\pi^\star=\pi^\star\partial\bar\partial$ and $\pi^\star$ is injective on forms, so $\partial\bar\partial\omega=0$. The second statement is proved in the same way with $\omega^{n-1}$ in place of $\omega$.
\end{proof}

\begin{Lem}[implications between the notions]\label{lem:impl}
Let $X$ be a compact complex manifold of dimension $n\geq2$. Then:
\begin{enumerate}[label=(\alph*)]
\item K\"ahler hyperbolic $\Rightarrow$ balanced hyperbolic;
\item balanced hyperbolic $\Rightarrow$ Gauduchon hyperbolic;
\item K\"ahler hyperbolic $\Rightarrow$ H-S hyperbolic;
\item H-S hyperbolic $\Rightarrow$ SKT hyperbolic.
\end{enumerate}
In each case the same metric $\omega$ works.
\end{Lem}

\begin{proof}
(a) Let $\omega$ be K\"ahler with $\pi^\star\omega=d\eta$, $\eta$ bounded. A K\"ahler metric is balanced. Since $d\eta=\pi^\star\omega$ is bounded and $d(d\eta)=0$,
\begin{eqnarray*}
\pi^\star\omega^{n-1} & = & (d\eta)^{n-1}\;=\;d\big(\eta\wedge(d\eta)^{n-2}\big),
\end{eqnarray*}
and $\eta\wedge(d\eta)^{n-2}$ is bounded. So $\omega^{n-1}$ is $\widetilde d$-bounded.

(b) Let $\omega$ be balanced with $\pi^\star\omega^{n-1}=d\Gamma$, $\Gamma$ a bounded $(2n-3)$-form. Decompose $\Gamma=\sum_{p+q=2n-3}\Gamma^{p,q}$; each component is bounded. Then
\begin{eqnarray*}
d\Gamma & = & \sum_{p+q=2n-3}\big(\partial\Gamma^{p,q}+\bar\partial\Gamma^{p,q}\big),
\end{eqnarray*}
and the component of bidegree $(n-1,n-1)$ comes from $\partial\Gamma^{n-2,n-1}$ and $\bar\partial\Gamma^{n-1,n-2}$ only. Since $\pi^\star\omega^{n-1}$ is of pure type $(n-1,n-1)$,
\begin{eqnarray*}
\pi^\star\omega^{n-1} & = & \partial\Gamma^{n-2,n-1}+\bar\partial\Gamma^{n-1,n-2},
\end{eqnarray*}
so $\omega^{n-1}$ is $\widetilde{(\partial+\bar\partial)}$-bounded; by Lemma~\ref{lem:auto}, $\omega$ is Gauduchon.

(c) Take $\rho=0$: then $\Omega=\omega$.

(d) Let $\Omega=\rho+\omega+\bar\rho$ with $\pi^\star\Omega=d\gamma$, $\gamma$ bounded, and write $\gamma=\gamma^{1,0}+\gamma^{0,1}$ (both bounded). The component of bidegree $(1,1)$ of $d\gamma$ is $\partial\gamma^{0,1}+\bar\partial\gamma^{1,0}$, and the component of bidegree $(1,1)$ of $\pi^\star\Omega$ is $\pi^\star\omega$. Hence
\begin{eqnarray*}
\pi^\star\omega & = & \partial\gamma^{0,1}+\bar\partial\gamma^{1,0}
\end{eqnarray*}
(one must take this component: $\Omega$ and $\omega$ are different forms). By Lemma~\ref{lem:auto}, $\omega$ is SKT, hence SKT hyperbolic.
\end{proof}

\begin{Rem}
By (d) and [Mar24, Thm.~3.5], H-S hyperbolic manifolds are Kobayashi hyperbolic, a result of [Khe26].
\end{Rem}

\subsection{Divisorial hyperbolicity}

A holomorphic map $f:\C^{n-1}\to X$ is {\bf non-degenerate at $x$} if $d_xf$ is injective (rank $n-1$). The set $\Sigma$ where $df$ has rank $<n-1$ is then a proper analytic subset, and $f^\star\omega\geq0$ is a Hermitian metric on $\C^{n-1}\setminus\Sigma$ and degenerates on $\Sigma$. Put $\omega_{n-1}=\omega^{n-1}/(n-1)!$, $\tau(z)=|z|^2$, $B_t=\{|z|<t\}$, $S_t=\{|z|=t\}$, and
\begin{eqnarray*}
\mathrm{Vol}_{\omega,f}(B_t) & = & \int_{B_t}f^\star\omega_{n-1},\qquad F(b)\;=\;\int_0^b\mathrm{Vol}_{\omega,f}(B_t)\,dt .
\end{eqnarray*}
The {\bf area measure} $d\sigma_{\omega,f,t}$ is the $(2n-3)$-dimensional Riemannian measure induced on $S_t\setminus\Sigma$ by $g=f^\star\omega$, and $\|d\tau\|_{f^\star\omega}$ is the norm of $d\tau$ for $g$. Write $A_{\omega,f}(S_t)=\int_{S_t}d\sigma_{\omega,f,t}$.

\begin{Def}[{[MP22]}]\label{def:growth}
\begin{enumerate}[label=(\arabic*)]
\item $f$ has {\bf subexponential growth} if
\begin{enumerate}[label=(\roman*)]
\item there are $C_1,r_0>0$ with
\begin{eqnarray}
\int_{S_t}\|d\tau\|_{f^\star\omega}\,d\sigma_{\omega,f,t} & \leq & C_1\,t\,\mathrm{Vol}_{\omega,f}(B_t)\qquad(t>r_0);\label{eq:growth-i}
\end{eqnarray}
\item for every $C>0$,
\begin{eqnarray}
\limsup_{b\to+\infty}\Big(\frac bC-\log F(b)\Big) & = & +\infty .\label{eq:growth-ii}
\end{eqnarray}
\end{enumerate}
\item $X$ is {\bf divisorially hyperbolic} if there is no holomorphic $f:\C^{n-1}\to X$, non-degenerate at some point, with subexponential growth.
\end{enumerate}
\end{Def}

Condition (i) compares the sphere $S_t$ and the ball $B_t$ as measured by $f^\star\omega$: for the Euclidean metric both sides are of order $t^{2n-1}$ and it holds. Condition (ii) says that $\log F(b)$ is, along a sequence, much smaller than $b/C$ for every $C$, i.e.\ that the volume of balls does not grow exponentially along that sequence.

\begin{Lem}[the two measures on the sphere]\label{lem:measures}
Let $f:\C^{n-1}\to X$ be non-degenerate at some point, $g=f^\star\omega$, and let $\nu=d\tau/\|d\tau\|_g$ on $\C^{n-1}\setminus\Sigma$. Then $d\sigma_{\omega,f,t}=\star_g\nu|_{S_t}$ and
\begin{eqnarray}
f^\star\omega_{n-1} & = & \nu\wedge\star_g\nu,\qquad \mathrm{Vol}_{\omega,f}(B_r)\;=\;\int_0^r\Big(\int_{S_t}d\mu_{\omega,f,t}\Big)dt,\qquad d\mu_{\omega,f,t}\;=\;\frac{2t}{\|d\tau\|_g}\,d\sigma_{\omega,f,t}.\label{eq:two-measures}
\end{eqnarray}
In particular the measure $d\mu_{\omega,f,t}$ that disintegrates the volume form $f^\star\omega_{n-1}$ on the spheres is {\bf not} the area measure $d\sigma_{\omega,f,t}$ induced by the degenerate metric $g$: they differ by the non-constant factor $2t/\|d\tau\|_g$.
\end{Lem}

\begin{proof}
For a Hermitian metric, $f^\star\omega_{n-1}=g^{n-1}/(n-1)!$ is the Riemannian volume form $dV_g$, and $\nu\wedge\star_g\nu=|\nu|_g^2\,dV_g=dV_g$. The $(2n-3)$-form $\star_g\nu$ restricted to the hypersurface $S_t$ with unit conormal $\nu$ is its Riemannian area form, i.e.\ $d\sigma_{\omega,f,t}$. Hence $dV_g=\frac{d\tau}{\|d\tau\|_g}\wedge d\sigma$, and since $d\tau=2t\,dt$ the coarea formula for $t=|z|$ gives \eqref{eq:two-measures}.
\end{proof}

\begin{Lem}[equivalent form of (ii)]\label{lem:liminf}
Let $f:\C^{n-1}\to X$ be non-degenerate at some point. Then (ii) holds if and only if
\begin{eqnarray*}
\liminf_{b\to+\infty}\frac{\log F(b)}{b} & = & 0 .
\end{eqnarray*}
\end{Lem}

\begin{proof}
$F$ is non-decreasing and positive for large $b$. If $F$ is bounded, $\log F(b)/b\to0$ and (ii) holds trivially. Otherwise $F\to\infty$, so $\log F\geq0$ for large $b$ and $\liminf\geq0$. If $\liminf=0$, take $b_j\to\infty$ with $\log F(b_j)\leq\varepsilon_jb_j$, $\varepsilon_j\to0$; then $\frac{b_j}C-\log F(b_j)\geq b_j(\frac1C-\varepsilon_j)\to+\infty$ for every $C$. If $\liminf=\ell>0$ (possibly $\infty$), choose $C>0$ with $1/C<\ell/2$ (or $1/C<1$ if $\ell=\infty$); then $\log F(b)\geq\frac\ell2b$ for large $b$ and $\frac bC-\log F(b)\to-\infty$, so (ii) fails.
\end{proof}

$X$ is {\bf Brody hyperbolic} if every holomorphic map $\C\to X$ is constant; for compact $X$ this is equivalent to Kobayashi hyperbolicity (Brody's theorem). Brody hyperbolic implies divisorially hyperbolic: a map $f:\C^{n-1}\to X$ non-degenerate at $x$ is non-constant on every complex line through $x$.

\begin{The}[{[MP22, Thm.~2.8], [Mar24, Thms.~3.5, 3.11]}]\label{thm:known}
Balanced hyperbolic $\Rightarrow$ divisorially hyperbolic. SKT hyperbolic $\Rightarrow$ Kobayashi hyperbolic. Gauduchon hyperbolic $\Rightarrow$ divisorially hyperbolic.
\end{The}

\section{Divisorial hyperbolicity of surfaces}\label{sec:surfaces}\label{sec:brody}

In this section $n=2$. Then $\C^{n-1}=\C$, $\omega_{n-1}=\omega$ and ``non-degenerate at some point'' means ``non-constant''. For $f:\C\to X$ non-constant let $\mu=f^\star\omega\geq0$ (a measure on $\C$, absolutely continuous) and
\begin{eqnarray*}
A(s) & = & \mu(B_s),\qquad F(b)\;=\;\int_0^bA(s)\,ds .
\end{eqnarray*}
The function $A$ is non-decreasing and positive for large $s$.

\begin{Prop}\label{prop:brody}
Let $X$ be a compact complex surface. Then $X$ is Brody hyperbolic if and only if $X$ is divisorially hyperbolic.
\end{Prop}

The implication ($\Rightarrow$) was recalled above. For ($\Leftarrow$) we show, by contraposition, that if $X$ is not Brody hyperbolic there is a non-constant $G:\C\to X$ of subexponential growth. The proof uses four lemmas.

\begin{Lem}[Condition (i) in dimension two]\label{lem:i}
For every non-constant $f:\C\to X$ and every $t>0$,
\begin{eqnarray}
\int_{S_t}\|d\tau\|_{f^\star\omega}\,d\sigma_{\omega,f,t} & = & 4\pi t^2 .\label{eq:4pit2}
\end{eqnarray}
Consequently (i) holds if and only if there are $c,r_0>0$ with $A(t)\geq c\,t$ for all $t>r_0$.
\end{Lem}

\begin{proof}
Write $f^\star\omega=m(z)\,\frac i2dz\wedge d\bar z$ with $m\geq0$ smooth, vanishing only on the discrete set $\Sigma$. Off $\Sigma$ the metric $g=m\,|dz|^2$ is conformal to the Euclidean one, so for $z\in S_t$
\begin{eqnarray*}
\|d\tau\|_g & = & \frac{|d\tau|_{\rm eucl}}{\sqrt m}\;=\;\frac{2t}{\sqrt m},\qquad d\sigma_{\omega,f,t}\;=\;\sqrt m\;ds ,
\end{eqnarray*}
where $ds$ is the Euclidean arc length (in real dimension $2$ the induced measure on the circle is the $g$-length element). The product is $2t\,ds$ and $S_t\cap\Sigma$ is finite, so the integral is $2t\cdot2\pi t=4\pi t^2$. Inserting \eqref{eq:4pit2} in \eqref{eq:growth-i} gives $4\pi t^2\leq C_1\,t\,A(t)$, i.e.\ $A(t)\geq(4\pi/C_1)\,t$.
\end{proof}

\begin{Lem}[Condition (ii)]\label{lem:ii}
If $F(b)>0$ for large $b$ and $\liminf_{b\to\infty}\frac{\log F(b)}{b}\leq0$, then \eqref{eq:growth-ii} holds.
\end{Lem}

\begin{proof}
There are $b_j\to\infty$ and $\varepsilon_j\to0$ with $\log F(b_j)\leq\varepsilon_jb_j$. For every $C>0$,
\begin{eqnarray*}
\frac{b_j}{C}-\log F(b_j) & \geq & b_j\Big(\frac1C-\varepsilon_j\Big)\;\longrightarrow\;+\infty .
\end{eqnarray*}
\end{proof}

\begin{Lem}[Precomposition by $z\mapsto z^k$]\label{lem:power}
Let $f:\C\to X$ be non-constant and $G=f\circ(z\mapsto z^k)$, $k\geq1$. Then $G$ is non-constant and
\begin{eqnarray*}
A_G(t) & = & k\,A(t^k)\qquad(t>0).
\end{eqnarray*}
\end{Lem}

\begin{proof}
The map $z\mapsto z^k$ is open and sends $B_t$ onto $B_{t^k}$ with multiplicity $k$, so
\begin{eqnarray*}
A_G(t) & = & \int_{B_t}(z^k)^\star(f^\star\omega)\;=\;k\int_{B_{t^k}}f^\star\omega\;=\;k\,A(t^k).
\end{eqnarray*}
\end{proof}

\begin{Lem}[Precomposition by $\exp$]\label{lem:exp}
Let $f:\C\to X$ be non-constant, $G=f\circ\exp$ and $N_t(w)=\#\{z\in B_t:\ e^z=w\}$ for $w\in\C^\ast$. Then $G$ is non-constant and
\begin{eqnarray}
A_G(t) & = & \int_{\C^\ast}N_t(w)\,d\mu(w).\label{eq:AG-exp}
\end{eqnarray}
Moreover:
\begin{enumerate}[label=(\alph*)]
\item {\bf (lower bound)} there are $a>0$ and $t_0$ with $A_G(t)\geq\frac a{2\pi}\,t$ for $t\geq t_0$;
\item {\bf (upper bound)} for every $t>0$,
\begin{eqnarray}
A_G(t) & \leq & \Big(1+\frac t\pi\Big)A(e^t).\label{eq:AG-upper}
\end{eqnarray}
\end{enumerate}
\end{Lem}

\begin{proof}
Formula \eqref{eq:AG-exp} is the change of variables $w=e^z$ with multiplicity.

(a) Since $\mu$ is absolutely continuous and $\mu(\C)>0$, there is a closed disc $D\subset\C^\ast$ with $a:=\mu(D)>0$. Put $R=\max_{w\in D}|\log|w||$. For $w=|w|e^{i\theta}\in D$ the solutions of $e^z=w$ are
\begin{eqnarray*}
z_m & = & \log|w|+i(\theta+2\pi m),\qquad m\in\Z,
\end{eqnarray*}
and $z_m\in B_t$ as soon as $|\theta+2\pi m|<\sqrt{t^2-R^2}$ (for $t>R$). An interval of length $2\sqrt{t^2-R^2}$ contains at least $\sqrt{t^2-R^2}/\pi-1$ points of $\theta+2\pi\Z$, hence
\begin{eqnarray*}
N_t(w) & \geq & \frac{\sqrt{t^2-R^2}}{\pi}-1\;\geq\;\frac t{2\pi}\qquad(t\ \mbox{large},\ w\in D).
\end{eqnarray*}
Then \eqref{eq:AG-exp} gives $A_G(t)\geq\int_DN_t\,d\mu\geq\frac{a}{2\pi}\,t$.

(b) If $z\in B_t$ and $e^z=w$, then $|{\rm Re}\,z|<t$ and $|{\rm Im}\,z|<t$, so $e^{-t}<|w|<e^t$, and the imaginary parts of the solutions form a progression of step $2\pi$ in $(-t,t)$, whence $N_t(w)\leq1+t/\pi$. Thus $N_t$ vanishes outside $\{e^{-t}<|w|<e^t\}$ and
\begin{eqnarray*}
A_G(t) & \leq & \Big(1+\frac t\pi\Big)\mu\big(\{e^{-t}<|w|<e^t\}\big)\;\leq\;\Big(1+\frac t\pi\Big)A(e^t).
\end{eqnarray*}
\end{proof}

\begin{proof}[Proof of Proposition~\ref{prop:brody}, ($\Leftarrow$)]
Assume $X$ is not Brody hyperbolic. By Brody's reparametrisation lemma ($X$ is compact) there is a non-constant $f:\C\to X$ with $\|f'\|_\omega\leq1$; hence $A(s)\leq C\,s^2$ for some $C>0$. Put
\begin{eqnarray*}
\lambda & = & \liminf_{s\to\infty}\frac{\log A(s)}{\log s}\;\in\;[0,2].
\end{eqnarray*}

{\it Case 1: $\lambda>0$.} Choose an integer $k\geq2/\lambda$ and $G=f\circ(z\mapsto z^k)$. For $s$ large $A(s)\geq s^{\lambda/2}$, so by Lemma~\ref{lem:power}, for $t\geq1$ large,
\begin{eqnarray*}
A_G(t) & = & k\,A(t^k)\;\geq\;k\,t^{k\lambda/2}\;\geq\;t ,
\end{eqnarray*}
which is (i) by Lemma~\ref{lem:i}. Also $A_G(t)\leq kC\,t^{2k}$, so $F_G(b)\leq kC\,b^{2k+1}$ and $\log F_G(b)=O(\log b)=o(b)$; since $F_G(b)\geq b^2/2$, Lemma~\ref{lem:ii} gives (ii).

{\it Case 2: $\lambda=0$.} Take $G=f\circ\exp$. By Lemma~\ref{lem:exp}(a) and Lemma~\ref{lem:i}, (i) holds (with $C_1=8\pi^2/a$). By \eqref{eq:AG-upper} and monotonicity of $A_G$,
\begin{eqnarray*}
F_G(b) & \leq & b\,A_G(b)\;\leq\;b\Big(1+\frac b\pi\Big)A(e^b).
\end{eqnarray*}
Since $\lambda=0$ there are $s_j\to\infty$, $\varepsilon_j\to0$ with $\log A(s_j)\leq\varepsilon_j\log s_j$. For $b_j=\log s_j$,
\begin{eqnarray*}
\log F_G(b_j) & \leq & \log b_j+\log\Big(1+\frac{b_j}\pi\Big)+\varepsilon_jb_j\;=\;o(b_j).
\end{eqnarray*}
Also $F_G(b)\geq c\,b^2>0$ by Lemma~\ref{lem:exp}(a), so Lemma~\ref{lem:ii} gives (ii).

In both cases $G:\C\to X$ is non-constant (Lemmas~\ref{lem:power} and~\ref{lem:exp}) and has subexponential growth, so $X$ is not divisorially hyperbolic.
\end{proof}

\begin{Rem}
Case 2 contains the case of bounded $A$ (e.g.\ $f:\C\to\mathbb P^1\subset X$ a rational curve). Brody's lemma only controls the order of $f$, i.e.\ condition (ii); condition (i) amounts, for $n=2$, to a linear growth of $A$ (Lemma~\ref{lem:i}) and is repaired by precomposition with $z^k$ or with $\exp$. This argument does not extend to $n\geq3$, where no Brody lemma is available for maps $\C^{n-1}\to X$ and the second term of the identity (8) of [MP22, Lemma~2.5] does not vanish.
\end{Rem}

\begin{Rem}
For $n=2$, [Mar24, Thm.~3.11] (Gauduchon hyperbolic $\Rightarrow$ divisorially hyperbolic) implies [Mar24, Thm.~3.5] (SKT hyperbolic $\Rightarrow$ Kobayashi hyperbolic): an SKT hyperbolic surface is Gauduchon hyperbolic, hence divisorially hyperbolic, hence Brody hyperbolic by Proposition~\ref{prop:brody}.
\end{Rem}

\section{Hermitian-symplectic hyperbolicity and divisorially K\texorpdfstring{\"a}{a}hler classes of surfaces}\label{sec:classes}

\subsection{The notions}

Let $X$ be a compact complex surface. We recall the notions of [MP22, \S3] in dimension two.

{\it De Rham classes and their Hodge components.} If $X$ is K\"ahler, every class $a\in H^2_{DR}(X,\mathbb R)$ decomposes uniquely as $a=a^{2,0}+a^{1,1}+a^{0,2}$ with $a^{0,2}=\overline{a^{2,0}}$ (Hodge decomposition), the components being represented by the corresponding components of the harmonic representative. We say $a$ has {\bf type $(1,1)$} if $a^{2,0}=0$; this is automatic if $h^{2,0}(X)=0$. The {\bf K\"ahler cone} $\mathcal K_X\subset H^{1,1}(X,\mathbb R)$ is the set of classes of K\"ahler forms.

{\it Aeppli cohomology.} $H^{1,1}_A(X,\mathbb R)$ is the space of real $(1,1)$-forms $\beta$ with $\partial\bar\partial\beta=0$, modulo the forms $\partial u+\bar\partial\bar u$ with $u$ a $(0,1)$-form. If $X$ is K\"ahler, the $\partial\bar\partial$-lemma gives $H^{1,1}_A(X,\mathbb R)\simeq H^{1,1}(X,\mathbb R)$.

{\it Gauduchon cone.} For $n=2$ a Gauduchon metric is a Hermitian metric $\omega$ with $\partial\bar\partial\omega=0$; its class $\{\omega\}_A\in H^{1,1}_A(X,\mathbb R)$ is well defined. The {\bf Gauduchon cone} $G_X\subset H^{1,1}_A(X,\mathbb R)$ is the set of such classes (an open convex cone).

\begin{Lem}[the map $P$; {[MP22, Lemma~3.1]} for $n=2$]\label{lem:P}
Let $\alpha$ be a real $d$-closed $2$-form. Then $\partial\bar\partial\alpha^{1,1}=0$, and replacing $\alpha$ by $\alpha+d\kappa$ ($\kappa$ a real $1$-form) changes $\alpha^{1,1}$ by $\partial\kappa^{0,1}+\bar\partial\kappa^{1,0}$. Hence
\begin{eqnarray*}
P:H^2_{DR}(X,\mathbb R)\;\longrightarrow\;H^{1,1}_A(X,\mathbb R),\qquad P(\{\alpha\}_{DR}) & = & \{\alpha^{1,1}\}_A
\end{eqnarray*}
is a well-defined linear map.
\end{Lem}

\begin{proof}
The $(1,2)$-component of $d\alpha=0$ is $\partial\alpha^{0,2}+\bar\partial\alpha^{1,1}=0$, so $\partial\bar\partial\alpha^{1,1}=-\partial\partial\alpha^{0,2}=0$. The $(1,1)$-component of $d\kappa=\partial\kappa^{1,0}+\bar\partial\kappa^{1,0}+\partial\kappa^{0,1}+\bar\partial\kappa^{0,1}$ is $\partial\kappa^{0,1}+\bar\partial\kappa^{1,0}$, which is zero in $H^{1,1}_A$ (as $\kappa^{1,0}=\overline{\kappa^{0,1}}$).
\end{proof}

\begin{Def}[{[MP22, Def.~3.6]}]
The {\bf divisorially K\"ahler cone} and the {\bf divisorially nef cone} are
\begin{eqnarray*}
DK_X & = & P^{-1}(G_X),\qquad DN_X\;=\;P^{-1}\big(\overline{G_X}\big)\;\subset\;H^2_{DR}(X,\mathbb R).
\end{eqnarray*}
\end{Def}

Finally let $V_X\subset H^2_{DR}(X,\mathbb R)$ be the subspace of classes having a $\widetilde d$-bounded representative. It is a linear subspace, and if one representative of a class is $\widetilde d$-bounded then all are, since $\pi^\star d\kappa=d(\pi^\star\kappa)$ with $\pi^\star\kappa$ bounded.

\subsection{An H-S form determines a K\texorpdfstring{\"a}{a}hler class}

\begin{Lem}[Hodge decomposition of an H-S form]\label{lem:decomp}
Let $X$ be a compact K\"ahler surface, $\omega_K$ a K\"ahler form and $\Omega=\rho+\omega+\bar\rho$ a real $d$-closed form with $\omega>0$. Then there are a holomorphic $2$-form $h$, a real $d$-closed $(1,1)$-form $\Omega_h^{1,1}$ and a real $1$-form $\kappa$ such that
\begin{eqnarray}
\Omega & = & h+\Omega_h^{1,1}+\bar h+d\kappa,\label{eq:dec1}\\
\rho & = & \bar h+\eta,\qquad \eta\;=\;\bar\partial\kappa^{0,1},\label{eq:dec2}\\
\omega & = & \Omega_h^{1,1}+\partial\kappa^{0,1}+\bar\partial\kappa^{1,0}.\label{eq:dec3}
\end{eqnarray}
The class $c=[\Omega_h^{1,1}]$ is the $(1,1)$-Hodge component of $[\Omega]$.
\end{Lem}

\begin{proof}
Write $\Omega=\Omega_h+d\kappa$ with $\Omega_h$ the $\omega_K$-harmonic representative. On a K\"ahler manifold the components of a harmonic form are harmonic, so $\Omega_h=h+\Omega_h^{1,1}+\bar h$ with $h$ a harmonic $(2,0)$-form, hence holomorphic. Taking the components of type $(0,2)$ and $(1,1)$ of \eqref{eq:dec1} gives \eqref{eq:dec2} and \eqref{eq:dec3}.
\end{proof}

\begin{Lem}[three numerical properties]\label{lem:num}
In the situation of Lemma~\ref{lem:decomp}, with $c=[\Omega_h^{1,1}]$:
\begin{enumerate}[label=(\alph*)]
\item for every irreducible curve $Y\subset X$, $\int_Yc=\int_Y\omega>0$;
\item $\int_Xc\wedge\omega_K=\int_X\omega\wedge\omega_K>0$;
\item $\int_Xc^2\;\geq\;\int_X\omega^2>0$.
\end{enumerate}
\end{Lem}

\begin{proof}
(a) The $(2,0)$ and $(0,2)$-components restrict to $0$ on a complex curve and $\int_Yd\kappa=0$, so by \eqref{eq:dec1} and $\rho|_Y=0$,
\begin{eqnarray*}
\int_Y\Omega_h^{1,1} & = & \int_Y\Omega\;=\;\int_Y\omega .
\end{eqnarray*}
(b) $\omega_K$ is closed and, by type, $h\wedge\omega_K=\bar h\wedge\omega_K=\rho\wedge\omega_K=0$, so $\int c\wedge\omega_K=\int\Omega\wedge\omega_K=\int\omega\wedge\omega_K$.

(c) By type, $\Omega_h^2=(\Omega_h^{1,1})^2+2h\wedge\bar h$ and $\Omega^2=\omega^2+2\rho\wedge\bar\rho$, and $\int\Omega_h^2=\int\Omega^2$ since $\Omega-\Omega_h=d\kappa$. Hence
\begin{eqnarray}
\int c^2 & = & \int\omega^2+2\int\rho\wedge\bar\rho-2\int h\wedge\bar h.\label{eq:c2}
\end{eqnarray}
Now $\rho=\bar h+\eta$ and $\bar\rho=h+\bar\eta$, so
\begin{eqnarray*}
\int\rho\wedge\bar\rho & = & \int\bar h\wedge h+\int\bar h\wedge\bar\eta+\int\eta\wedge h+\int\eta\wedge\bar\eta .
\end{eqnarray*}
The two cross terms vanish: $\bar\eta=\partial\kappa^{1,0}$ and $\bar h\wedge\bar\partial\kappa^{1,0}=0$ by type, so $\bar h\wedge\bar\eta=d(\bar h\wedge\kappa^{1,0})$ because $d\bar h=0$; similarly $\eta\wedge h=d(\kappa^{0,1}\wedge h)$. Moreover $u\wedge\bar u\geq0$ for every $(2,0)$-form $u$, so $\int\eta\wedge\bar\eta\geq0$. Therefore $\int\rho\wedge\bar\rho\geq\int h\wedge\bar h$, and \eqref{eq:c2} gives (c).
\end{proof}

\begin{Lem}[K\"ahler class]\label{lem:Kclass}
In the situation of Lemma~\ref{lem:decomp}, $c$ is a K\"ahler class.
\end{Lem}

\begin{proof}
By Demailly--P\u{a}un [DP04], the K\"ahler cone of a compact K\"ahler manifold is a connected component of the set $\mathcal P$ of classes $\alpha$ with $\int_Y\alpha^p>0$ for every irreducible analytic subset $Y$ of dimension $p$. For $s\in[0,1]$ put $c_s=(1-s)[\omega_K]+s\,c$. For a curve $Y$, $\int_Yc_s=(1-s)\int_Y\omega_K+s\int_Y c>0$ by Lemma~\ref{lem:num}(a). For $Y=X$,
\begin{eqnarray*}
\int c_s^2 & = & (1-s)^2\!\int\omega_K^2+2s(1-s)\!\int c\wedge\omega_K+s^2\!\int c^2\;>\;0
\end{eqnarray*}
by Lemma~\ref{lem:num}(b),(c). So the whole segment lies in $\mathcal P$; it is connected and contains $[\omega_K]$, hence $c$ lies in the component of $\mathcal P$ containing $[\omega_K]$: $c$ is a K\"ahler class.
\end{proof}

\subsection{Explicit description of $DK_X$ and $DN_X$}

\begin{Prop}\label{prop:DK}
Let $X$ be a compact complex surface with $DK_X\neq\emptyset$. Then:
\begin{enumerate}[label=(\roman*)]
\item every class of $DK_X$ contains an H-S form;
\item $X$ is K\"ahler, $P(a)=a^{1,1}$ (via $H^{1,1}_A\simeq H^{1,1}(X,\mathbb R)$) and $G_X=\mathcal K_X$;
\item $DK_X=\{a:\ a^{1,1}\in\mathcal K_X\}$, $DN_X=\{a:\ a^{1,1}\in\overline{\mathcal K_X}\}$ and $\overline{DK_X}=DN_X$.
\end{enumerate}
\end{Prop}

\begin{proof}
(i) Let $a\in DK_X$ and $\Omega_0$ a closed real representative. Since $P(a)=\{\Omega_0^{1,1}\}_A\in G_X$ there is a Gauduchon metric $\omega$ and a $(0,1)$-form $u$ with
\begin{eqnarray*}
\omega & = & \Omega_0^{1,1}+\partial u+\bar\partial\bar u
\end{eqnarray*}
(everything is real, so the second potential may be taken to be $\bar u$). Put $\Omega=\Omega_0+d(u+\bar u)$: it is closed, cohomologous to $\Omega_0$, and by Lemma~\ref{lem:P} its $(1,1)$-component is $\omega>0$. So $\Omega$ is an H-S form.

(ii) By (i), $X$ carries a Hermitian-symplectic metric, hence is K\"ahler ([ST10]: a compact complex surface carrying an H-S metric is K\"ahler). On a K\"ahler surface, $P(a)=\{a_h^{1,1}\}_A$ for the harmonic representative, i.e.\ $P(a)=a^{1,1}$ in $H^{1,1}(X,\mathbb R)$. A K\"ahler form is Gauduchon, so $\mathcal K_X\subset G_X$. Conversely let $g\in G_X$, viewed as a class of $H^{1,1}(X,\mathbb R)\subset H^2_{DR}(X,\mathbb R)$. Then $P(g)=g\in G_X$, so $g\in DK_X$, $g$ contains an H-S form by (i), and by Lemma~\ref{lem:Kclass} its $(1,1)$-component $g$ is K\"ahler.

(iii) The first two identities follow from (ii), since $G_X=\mathcal K_X$ is open. The inclusion $\overline{DK_X}\subset DN_X$ holds because $DN_X=P^{-1}(\overline{G_X})$ is closed ($P$ is continuous) and contains $DK_X$. Conversely let $a^{1,1}\in\overline{\mathcal K_X}$ and choose $k_j\in\mathcal K_X$ with $k_j\to a^{1,1}$. Put
\begin{eqnarray*}
a_j & = & a+(k_j-a^{1,1}) .
\end{eqnarray*}
Then $a_j\to a$ and $a_j^{1,1}=k_j\in\mathcal K_X$, so $a_j\in DK_X$ and $a\in\overline{DK_X}$.
\end{proof}

\begin{Cor}[Question~3.9 of {[MP22]} in dimension two]\label{cor:Q39}
If $X$ is a compact complex surface with $DK_X\neq\emptyset$, then $\overline{DK_X}=DN_X$.
\end{Cor}

\begin{Question}
For $n\geq3$, is $\overline{DK_X}=DN_X$ when $DK_X\neq\emptyset$ (Question~3.9 of [MP22])?
\end{Question}

\subsection{Hermitian-symplectic hyperbolicity}

\begin{The}\label{thm:HSDK}
Let $X$ be a compact complex surface.
\begin{enumerate}[label=(\roman*)]
\item $X$ is H-S hyperbolic if and only if $V_X\cap DK_X\neq\emptyset$.
\item $X$ is K\"ahler hyperbolic if and only if $V_X\cap\mathcal K_X\neq\emptyset$.
\end{enumerate}
When $X$ is K\"ahler, $\mathcal K_X=DK_X\cap H^{1,1}(X,\mathbb R)$.
\end{The}

\begin{proof}
(i) If $\Omega$ is an H-S hyperbolic form with $(1,1)$-part $\omega$, then $[\Omega]\in V_X$ and $P([\Omega])=\{\omega\}_A\in G_X$ ($\omega$ is Gauduchon by Lemma~\ref{lem:auto}); so $[\Omega]\in V_X\cap DK_X$. Conversely, if $a\in V_X\cap DK_X$, Proposition~\ref{prop:DK}(i) gives an H-S form $\Omega\in a$, which is $\widetilde d$-bounded since $a\in V_X$; thus $\Omega^{1,1}$ is an H-S hyperbolic metric. (ii) is the definition of K\"ahler hyperbolicity, and the last statement is Proposition~\ref{prop:DK}(iii).
\end{proof}

\begin{Cor}\label{cor:K}
Let $X$ be an H-S hyperbolic compact complex surface. If $V_X\cap DK_X$ contains a class of type $(1,1)$ (this holds as soon as $h^{2,0}(X)=0$), then $X$ is K\"ahler hyperbolic. Hence, for surfaces with $h^{2,0}=0$,
\begin{eqnarray*}
\mbox{H-S hyperbolic} & \iff & \mbox{K\"ahler hyperbolic}.
\end{eqnarray*}
\end{Cor}

\begin{proof}
Let $a\in V_X\cap DK_X$ have type $(1,1)$. Then $a=a^{1,1}\in\mathcal K_X$ by Proposition~\ref{prop:DK}(iii), so $V_X\cap\mathcal K_X\neq\emptyset$ and Theorem~\ref{thm:HSDK}(ii) applies.
\end{proof}

\begin{Rem}\label{rem:missing}
For $a\in V_X\cap DK_X$ the class $a^{1,1}$ is always K\"ahler, and $X$ is K\"ahler hyperbolic if and only if $a^{1,1}\in V_X$, i.e.\ if and only if the real class $a^{2,0}+a^{0,2}$ (the class of $h+\bar h$) is $\widetilde d$-bounded. This is the only missing point when $h^{2,0}(X)>0$.
\end{Rem}

\begin{Question}\label{q:HS}
Is every H-S hyperbolic compact complex manifold K\"ahler hyperbolic? (Known for surfaces with $h^{2,0}=0$, Corollary~\ref{cor:K}.)
\end{Question}

\subsection{The canonical bundle}

\begin{Prop}[Question~1.5 of {[MP22]} for K\"ahler surfaces]\label{prop:Q15}
Let $X$ be a compact K\"ahler surface that is divisorially hyperbolic (equivalently Brody hyperbolic). Then $K_X$ is nef; in particular $c_1(K_X)\in DN_X$, i.e.\ $K_X$ is divisorially nef.
\end{Prop}

\begin{proof}
By Proposition~\ref{prop:brody} $X$ contains no entire curve, hence no rational curve; so $X$ has no $(-1)$-curve and is minimal. A K\"ahler surface of Kodaira dimension $-\infty$ is $\mathbb P^2$ or ruled, hence contains rational curves; so $\kappa(X)\geq0$, and the canonical class of a minimal K\"ahler surface with $\kappa\geq0$ is nef (see [BHPV04, Chap.~VII--VIII]). The class $c_1(K_X)$ has type $(1,1)$, so $P(c_1(K_X))=c_1(K_X)\in\overline{\mathcal K_X}$, i.e.\ $c_1(K_X)\in DN_X$ by Proposition~\ref{prop:DK}(iii).
\end{proof}

\begin{Question}
For $n\geq3$, if $X$ is balanced hyperbolic, is $K_X$ divisorially nef (Question~1.5 of [MP22])?
\end{Question}

\section{Codimension one in arbitrary dimension}\label{sec:codim}

\subsection{Compact hypersurfaces that lift to the universal cover}

\begin{Prop}\label{prop:hyper}
Let $X$ be a compact complex manifold of dimension $n\geq2$ and $\omega$ a Hermitian metric such that $\pi^\star\omega^{n-1}=\partial\Gamma+\bar\partial\Gamma'$ on $\widetilde X$ for some smooth forms $\Gamma,\Gamma'$ (this holds if $\omega$ is Gauduchon hyperbolic, hence if it is balanced hyperbolic; no boundedness is needed). Let $Y$ be a compact connected complex manifold of dimension $n-1$ and $f:Y\to X$ a holomorphic map with $f_\ast\pi_1(Y)=\{1\}$. Then $df$ has rank $<n-1$ at every point of $Y$; in particular $f$ is not generically finite onto its image.
\end{Prop}

\begin{proof}
Since $f_\ast\pi_1(Y)$ is trivial, $f=\pi\circ\tilde f$ for a holomorphic $\tilde f:Y\to\widetilde X$ (lifting criterion). The form $\tilde f^\star\Gamma$ has bidegree $(n-2,n-1)$ on the $(n-1)$-dimensional manifold $Y$, so $\bar\partial\tilde f^\star\Gamma$ has bidegree $(n-2,n)$ and vanishes; hence $\tilde f^\star\partial\Gamma=\partial\tilde f^\star\Gamma=d\tilde f^\star\Gamma$, and likewise for $\Gamma'$. By Stokes ($Y$ is compact without boundary),
\begin{eqnarray*}
\int_Yf^\star\omega^{n-1} & = & \int_Y\tilde f^\star\big(\pi^\star\omega^{n-1}\big)\;=\;\int_Yd\big(\tilde f^\star\Gamma+\tilde f^\star\Gamma'\big)\;=\;0 .
\end{eqnarray*}
But $f^\star\omega^{n-1}\geq0$ pointwise and is $>0$ exactly where $df$ has rank $n-1$; so it vanishes identically.
\end{proof}

\subsection{A linear isoperimetric inequality on the universal cover}

Recall that a group $\Gamma$ is {\bf amenable} if it carries a finitely additive $\Gamma$-invariant probability measure on all its subsets. Finite groups, abelian, nilpotent and solvable groups are amenable; free groups of rank $\geq2$ and fundamental groups of compact Riemann surfaces of genus $\geq2$ are not. For a complete Riemannian manifold $M$ the {\bf Cheeger constant} is
\begin{eqnarray*}
h(M) & = & \inf_\Omega\frac{\mathrm{area}(\partial\Omega)}{\mathrm{vol}(\Omega)},
\end{eqnarray*}
the infimum over bounded open sets $\Omega$ with $C^1$ boundary; $h(M)>0$ is a {\bf linear isoperimetric inequality} $\mathrm{vol}(\Omega)\leq h(M)^{-1}\mathrm{area}(\partial\Omega)$. By Brooks' theorem [Bro81], for a compact manifold $X$, $h(\widetilde X)>0$ if and only if $\pi_1(X)$ is non-amenable.

\begin{Lem}\label{lem:iso}
Let $X$ be a compact K\"ahler manifold that is Gauduchon hyperbolic (in particular, balanced hyperbolic) and let $\omega_K$ be any K\"ahler metric on $X$. There is $C>0$ such that, for every bounded open $\Omega\subset\widetilde X$ with $C^1$ boundary,
\begin{eqnarray*}
\mathrm{vol}_{\widetilde\omega_K}(\Omega) & \leq & C\;\mathrm{area}_{\widetilde\omega_K}(\partial\Omega).
\end{eqnarray*}
\end{Lem}

\begin{proof}
Let $\omega$ be a Gauduchon hyperbolic metric: $\pi^\star\omega^{n-1}=\partial\Gamma+\bar\partial\Gamma'$ with $\Gamma,\Gamma'$ bounded. By bidegree $\bar\partial\Gamma=0$ and $\partial\Gamma'=0$, and $d\widetilde\omega_K=0$, so
\begin{eqnarray*}
d\big((\Gamma+\Gamma')\wedge\widetilde\omega_K\big) & = & \pi^\star\omega^{n-1}\wedge\widetilde\omega_K .
\end{eqnarray*}
Since $X$ is compact and both $\omega^{n-1}\wedge\omega_K$ and $\omega_K^n/n!$ are positive $(n,n)$-forms, there is $c>0$ with $\omega^{n-1}\wedge\omega_K\geq c\,\omega_K^n/n!$. By Stokes, with $C_0=\sup|(\Gamma+\Gamma')\wedge\widetilde\omega_K|$,
\begin{eqnarray*}
c\,\mathrm{vol}(\Omega) & \leq & \int_\Omega\pi^\star\omega^{n-1}\wedge\widetilde\omega_K\;=\;\int_{\partial\Omega}(\Gamma+\Gamma')\wedge\widetilde\omega_K\;\leq\;C_0\,\mathrm{area}(\partial\Omega).
\end{eqnarray*}
\end{proof}

\begin{Cor}\label{cor:amen}
A compact K\"ahler manifold that is balanced hyperbolic or Gauduchon hyperbolic has a non-amenable fundamental group; in particular $\pi_1(X)$ is infinite.
\end{Cor}

\begin{proof}
Lemma~\ref{lem:iso} gives $h(\widetilde X)>0$; apply Brooks' theorem.
\end{proof}

\subsection{Maps from \texorpdfstring{$\C^{n-1}$}{C^(n-1)}: isoperimetric inequality and exponential growth}

Let $\omega$ be Gauduchon hyperbolic: $\pi^\star\omega_{n-1}=\partial\Gamma+\bar\partial\Gamma'$ with $\Gamma,\Gamma'$ bounded, and put $\eta=\Gamma+\Gamma'$ (a bounded $(2n-3)$-form). Let $f:\C^{n-1}\to X$ be non-degenerate at some point and $\tilde f$ a lift to $\widetilde X$ ($\C^{n-1}$ is simply connected). As in the proof of Proposition~\ref{prop:hyper}, $\bar\partial\tilde f^\star\Gamma=0$ and $\partial\tilde f^\star\Gamma'=0$ by bidegree, so
\begin{eqnarray}
f^\star\omega_{n-1} & = & d\big(\tilde f^\star\eta\big)\quad\mbox{on }\C^{n-1}.\label{eq:fstar-exact}
\end{eqnarray}

\begin{Lem}[comparison with the area measure]\label{lem:comp}
For every $r>0$, with $C_\eta=\sup\|\eta\|_{\widetilde\omega}$,
\begin{eqnarray*}
\Big|\int_{S_r}\tilde f^\star\eta\Big| & \leq & C_\eta\,A_{\omega,f}(S_r).
\end{eqnarray*}
\end{Lem}

\begin{proof}
For tangent vectors $v_1,\dots,v_{2n-3}$ to $S_r$ at a point outside $\Sigma$, forming an orthonormal basis for $g=f^\star\omega=\tilde f^\star\widetilde\omega$, the vectors $\tilde f_\ast v_j$ are orthonormal in $(T\widetilde X,\widetilde\omega)$, so
\begin{eqnarray*}
\big|(\tilde f^\star\eta)(v_1,\dots,v_{2n-3})\big| & \leq & \|\eta\|_{\widetilde\omega}\,\big|\tilde f_\ast v_1\wedge\dots\wedge\tilde f_\ast v_{2n-3}\big|_{\widetilde\omega}\;=\;\|\eta\|_{\widetilde\omega}.
\end{eqnarray*}
Hence $|\tilde f^\star\eta|\leq C_\eta\,d\sigma_{\omega,f,r}$ on $S_r\setminus\Sigma$; and $\Sigma\cap S_r$ has $(2n-3)$-dimensional measure zero because $\Sigma$ has real dimension $\leq2n-4$.
\end{proof}

\begin{Prop}[isoperimetric inequality for maps]\label{prop:iso}
There is $C>0$ (depending only on $\omega$) such that, for every $f:\C^{n-1}\to X$ non-degenerate at some point and every $r>0$,
\begin{eqnarray*}
\mathrm{Vol}_{\omega,f}(B_r) & \leq & C\,A_{\omega,f}(S_r).
\end{eqnarray*}
\end{Prop}

\begin{proof}
By \eqref{eq:fstar-exact}, Stokes and Lemma~\ref{lem:comp},
\begin{eqnarray*}
\mathrm{Vol}_{\omega,f}(B_r) & = & \int_{B_r}d\big(\tilde f^\star\eta\big)\;=\;\int_{S_r}\tilde f^\star\eta\;\leq\;C_\eta\,A_{\omega,f}(S_r).
\end{eqnarray*}
\end{proof}

\begin{Lem}[coarea and Cauchy--Schwarz]\label{lem:coarea}
For every $f:\C^{n-1}\to X$ non-degenerate at some point and $r>0$,
\begin{eqnarray*}
\mathrm{Vol}_{\omega,f}(B_r) & \geq & 2\int_0^r\frac{A_{\omega,f}(S_t)^2}{\displaystyle\int_{S_t}\|d\tau\|\,d\sigma_{\omega,f,t}}\;t\,dt .
\end{eqnarray*}
\end{Lem}

\begin{proof}
The coarea formula for $\tau=|z|^2$ and the metric $g=f^\star\omega$ gives, with $d\tau=2t\,dt$,
\begin{eqnarray*}
\mathrm{Vol}_{\omega,f}(B_r) & = & \int_0^{r^2}\Big(\int_{S_{\sqrt\tau}}\frac{d\sigma}{\|d\tau\|}\Big)d\tau\;=\;2\int_0^r\Big(\int_{S_t}\frac{d\sigma_{\omega,f,t}}{\|d\tau\|}\Big)t\,dt .
\end{eqnarray*}
By the Cauchy--Schwarz inequality on $S_t$,
\begin{eqnarray*}
A_{\omega,f}(S_t)^2 & = & \Big(\int_{S_t}\|d\tau\|^{1/2}\,\|d\tau\|^{-1/2}\,d\sigma\Big)^2\;\leq\;\int_{S_t}\|d\tau\|\,d\sigma\;\cdot\;\int_{S_t}\frac{d\sigma}{\|d\tau\|},
\end{eqnarray*}
which gives the lower bound on the inner integral.
\end{proof}

\begin{Prop}[exponential growth]\label{prop:exp}
Let $\omega$ be Gauduchon hyperbolic and let $f:\C^{n-1}\to X$ be non-degenerate at some point and satisfy (i) with constants $C_1,r_0$ ($r_0$ larger than the norm of a non-degenerate point). Put
\begin{eqnarray*}
F_{r_0}(b) & = & \int_{r_0}^b\mathrm{Vol}_{\omega,f}(B_t)\,dt,\qquad c_0\;=\;\frac{2}{C_1C^2},
\end{eqnarray*}
with $C$ as in Proposition~\ref{prop:iso}. Then
\begin{eqnarray}
F_{r_0}(b) & \geq & F_{r_0}(a)\,e^{c_0(b-a)}\qquad(b>a>r_0).\label{eq:gronwall}
\end{eqnarray}
In particular $f$ does not have subexponential growth, which reproves [MP22, Thm.~2.8] and [Mar24, Thm.~3.11]. Note that the proof uses the two measures of Lemma~\ref{lem:measures} separately: the area measure $d\sigma$ in the Stokes bound (Proposition~\ref{prop:iso}) and the disintegration measure $d\mu=(2t/\|d\tau\|)\,d\sigma$ in the coarea formula (Lemma~\ref{lem:coarea}).
\end{Prop}

\begin{proof}
For $t>r_0$, Proposition~\ref{prop:iso} gives $A_{\omega,f}(S_t)^2\geq\mathrm{Vol}(B_t)^2/C^2$, and (i) gives $\int_{S_t}\|d\tau\|\,d\sigma\leq C_1t\,\mathrm{Vol}(B_t)$. Hence the integrand in Lemma~\ref{lem:coarea} satisfies
\begin{eqnarray*}
\frac{A(S_t)^2}{\int_{S_t}\|d\tau\|\,d\sigma}\,t & \geq & \frac{\mathrm{Vol}(B_t)^2/C^2}{C_1t\,\mathrm{Vol}(B_t)}\,t\;=\;\frac{\mathrm{Vol}(B_t)}{C_1C^2},
\end{eqnarray*}
and it is $\geq0$ for $t\leq r_0$. Therefore, for $r>r_0$,
\begin{eqnarray*}
F_{r_0}'(r)\;=\;\mathrm{Vol}(B_r) & \geq & \frac2{C_1C^2}\int_{r_0}^r\mathrm{Vol}(B_t)\,dt\;=\;c_0\,F_{r_0}(r),
\end{eqnarray*}
and Gronwall's lemma gives \eqref{eq:gronwall} ($F_{r_0}>0$ for $r>r_0$ since $f$ is non-degenerate somewhere). If $f$ also satisfied (ii), take $C'=2/c_0$ in \eqref{eq:growth-ii}: since $F\geq F_{r_0}$, for $b>a>r_0$
\begin{eqnarray*}
\frac b{C'}-\log F(b) & \leq & \frac{c_0}2\,b-\log F_{r_0}(a)-c_0(b-a)\;\longrightarrow\;-\infty,
\end{eqnarray*}
a contradiction.
\end{proof}

\section{Examples}\label{sec:ex}

\subsection{Projective spaces are not divisorially hyperbolic}

In [MP22, \S2.3 (V)] it is shown that the linear embedding $j:\C^{n-1}\to\mathbb P^n$, $z\mapsto[1:z:0]$, satisfies (ii) but that the available estimate falls short of (i), so the question was left open there. As for $n=2$ (Section~\ref{sec:brody}), precomposing with the exponential repairs (i). Put $m=n-1$,
\begin{eqnarray*}
E(z) & = & [1:e^{z_1}:\dots:e^{z_m}:0]\in\mathbb P^n,\qquad z_j=x_j+iy_j,
\end{eqnarray*}
$\omega=i\partial\bar\partial\log|Z|^2$ the Fubini--Study metric (the choice of metric is irrelevant by [MP22, Lemma~2.6]), and
\begin{eqnarray*}
u_j\;=\;e^{2x_j},\qquad S\;=\;1+\sum_{j=1}^mu_j,\qquad p_0\;=\;\frac1S,\qquad p_j\;=\;\frac{u_j}S,\qquad P\;=\;p_0p_1\cdots p_m .
\end{eqnarray*}

\begin{Lem}[the pull-back metric]\label{lem:Pn-metric}
$E^\star\omega=i\sum_{j,k}h_{j\bar k}\,dz_j\wedge d\bar z_k$ with
\begin{eqnarray}
h_{j\bar k} & = & p_j\delta_{jk}-p_jp_k,\qquad \det h\;=\;P,\qquad h^{-1}\;=\;\mathrm{diag}\Big(\frac1{p_j}\Big)+\frac1{p_0}\,\mathbf 1\mathbf 1^T,\label{eq:h}
\end{eqnarray}
so $E$ is non-degenerate at every point, $\mathrm{tr}\,h^{-1}=\sum_{j=1}^m\big(\frac1{p_j}+\frac1{p_0}\big)$ and
\begin{eqnarray}
E^\star\omega_m & = & 2^m\,P\,dx\,dy .\label{eq:volform}
\end{eqnarray}
\end{Lem}

\begin{proof}
$E^\star\omega=i\partial\bar\partial\varphi$ with $\varphi=\log S$ depending only on $x$, so $h_{j\bar k}=\frac14\partial_{x_j}\partial_{x_k}\varphi$. Since $\partial_{x_j}\varphi=2p_j$ and $\partial_{x_k}u_j=2u_j\delta_{jk}$,
\begin{eqnarray*}
\frac14\,\partial_{x_j}\partial_{x_k}\varphi & = & p_j\delta_{jk}-p_jp_k .
\end{eqnarray*}
Write $h=D-pp^T$ with $D=\mathrm{diag}(p_j)$. Then $\det h=\det D\,(1-p^TD^{-1}p)=\prod_jp_j\,(1-\sum_jp_j)=P$, and the Sherman--Morrison formula gives $h^{-1}=D^{-1}+D^{-1}pp^TD^{-1}/(1-\sum p_j)=D^{-1}+p_0^{-1}\mathbf 1\mathbf 1^T$. Finally $\omega_m=\det h\prod_j\big(i\,dz_j\wedge d\bar z_j\big)$ and $i\,dz\wedge d\bar z=2\,dx\wedge dy$.
\end{proof}

\begin{Lem}[volume and condition (ii)]\label{lem:Pn-vol}
Let $V_m$ be the volume of the unit ball of $\mathbb R^m$. Then
\begin{eqnarray}
\mathrm{Vol}_{\omega,E}(B_t) & = & 2^mV_m\int_{|x|<t}P(x)\,\big(t^2-|x|^2\big)^{m/2}dx,\qquad\int_{\mathbb R^m}P\,dx\;=\;\frac{2^{-m}}{m!}.\label{eq:Pn-vol}
\end{eqnarray}
Consequently $c\,t^m\leq\mathrm{Vol}_{\omega,E}(B_t)\leq C\,t^m$ for $t$ large, $F(b)\leq C\,b^{m+1}$, and $E$ satisfies (ii).
\end{Lem}

\begin{proof}
Integrating \eqref{eq:volform} over $|y|<\sqrt{t^2-|x|^2}$ gives the first identity. With $P=\prod u_j/S^{m+1}$ and $dx_j=du_j/(2u_j)$,
\begin{eqnarray*}
\int_{\mathbb R^m}P\,dx & = & 2^{-m}\int_{\mathbb R_+^m}\frac{du}{(1+\sum u_j)^{m+1}}\;=\;\frac{2^{-m}}{(m-1)!}\int_0^\infty\frac{r^{m-1}}{(1+r)^{m+1}}\,dr\;=\;\frac{2^{-m}}{m!}.
\end{eqnarray*}
The upper bound follows from $(t^2-|x|^2)^{m/2}\leq t^m$; for the lower bound restrict the integral to $|x|<t/\sqrt2$, where $(t^2-|x|^2)^{m/2}\geq(t^2/2)^{m/2}$ and $\int_{|x|<t/\sqrt2}P\,dx\geq\frac12\int P\,dx$ for $t$ large. Then $F(b)\leq Cb^{m+1}$, $\log F(b)=O(\log b)=o(b)$, and Lemma~\ref{lem:ii} applies.
\end{proof}

\begin{Lem}[the trace term]\label{lem:Pn-trace}
Let $g=P\,\mathrm{tr}\,h^{-1}=\sum_{j=1}^m\frac P{p_j}+m\frac P{p_0}$, where $\frac P{p_j}=\prod_{i\neq j}u_i/S^m$ ($1\leq j\leq m$) and $\frac P{p_0}=\prod_iu_i/S^m$. Then $\int_{|x|<t}g\,dx\leq C(1+t)$.
\end{Lem}

\begin{proof}
Fix $j\geq1$ and put $a=1+u_j$. Integrating first over the variables $x_i$, $i\neq j$, with $dx_i=du_i/(2u_i)$,
\begin{eqnarray*}
\int_{\mathbb R^{m-1}}\frac P{p_j}\prod_{i\neq j}dx_i & = & 2^{-(m-1)}\int_{\mathbb R_+^{m-1}}\frac{du}{(a+\sum_{i\neq j}u_i)^m}\;=\;\frac{2^{-(m-1)}}{(m-1)!\,(1+u_j)},
\end{eqnarray*}
and $\int_{-t}^t\frac{dx_j}{1+e^{2x_j}}\leq t+\frac12$. For $P/p_0$, since $u_i<e^{2t}$ on $B_t$, we have $\sum u_i<me^{2t}$ and
\begin{eqnarray*}
\int_{|x|<t}\frac P{p_0}\,dx & \leq & 2^{-m}\int_{\{\sum u_i<me^{2t}\}}\frac{du}{(1+\sum u_i)^m}\;=\;\frac{2^{-m}}{(m-1)!}\int_0^{me^{2t}}\frac{r^{m-1}}{(1+r)^m}\,dr\;\leq\;\frac{2^{-m}\log(1+me^{2t})}{(m-1)!}.
\end{eqnarray*}
Summing gives $\int_{|x|<t}g\,dx\leq C(1+t)$.
\end{proof}

\begin{Prop}\label{prop:Pn}
For every $n\geq2$ the projective space $\mathbb P^n$ is not divisorially hyperbolic: the map $E$ is non-degenerate and has subexponential growth.
\end{Prop}

\begin{proof}
By Lemma~\ref{lem:Pn-metric}, $E$ is non-degenerate everywhere; by Lemma~\ref{lem:Pn-vol}, (ii) holds. For (i): $E^\star\omega$ is K\"ahler, so the second term of formula (8) of [MP22, Lemma~2.5] vanishes and, $\Lambda$ denoting the trace for $g=E^\star\omega$,
\begin{eqnarray*}
\int_{S_t}\|d\tau\|\,d\sigma_{\omega,E,t} & = & 2\int_{B_t}\Lambda\big(i\partial\bar\partial\tau\big)\,E^\star\omega_m\;=\;2^{m+1}V_m\int_{|x|<t}g(x)\big(t^2-|x|^2\big)^{m/2}dx ,
\end{eqnarray*}
because $i\partial\bar\partial\tau=\sum i\,dz_j\wedge d\bar z_j$, $\Lambda(i\partial\bar\partial\tau)=\mathrm{tr}\,h^{-1}$ and \eqref{eq:volform}. By Lemma~\ref{lem:Pn-trace} the right-hand side is at most $C\,t^m(1+t)\leq C'\,t^{m+1}$, while $t\,\mathrm{Vol}_{\omega,E}(B_t)\geq c\,t^{m+1}$ by Lemma~\ref{lem:Pn-vol}. This is (i).
\end{proof}

\begin{Rem}
For $n=2$ ($m=1$), $g\equiv1$ and the integral is exactly $4\pi t^2$, in agreement with Lemma~\ref{lem:i}. Since $\mathbb P^{n-1}\subset\mathbb P^n$, $\mathbb P^n$ is also not balanced hyperbolic by [MP22, Prop.~2.10] (or by Proposition~\ref{prop:hyper}).
\end{Rem}

\subsection{Divisorially hyperbolic manifolds that are not balanced hyperbolic}

By Theorem~\ref{thm:known} balanced (and Gauduchon) hyperbolic manifolds are divisorially hyperbolic. We show that the converse fails in every dimension $n\geq2$. The obstructions to balanced and Gauduchon hyperbolicity are Corollary~\ref{cor:amen} (amenable fundamental group) and Proposition~\ref{prop:hyper} (liftable hypersurface).

\begin{Prop}\label{prop:criteria}
Let $X$ be a compact K\"ahler manifold that is Brody hyperbolic. Then $X$ is divisorially hyperbolic. If moreover $\pi_1(X)$ is amenable (e.g.\ finite or abelian), or if $X$ contains a liftable hypersurface (a map $f:Y\to X$ as in Proposition~\ref{prop:hyper}, generically finite), then $X$ is neither balanced hyperbolic nor Gauduchon hyperbolic.
\end{Prop}

\begin{proof}
The first claim was recalled in Section~\ref{sec:prelim}; the others are Corollary~\ref{cor:amen} and Proposition~\ref{prop:hyper}.
\end{proof}

\begin{Lem}[products]\label{lem:prod}
Let $X_1$ be a compact Brody hyperbolic manifold of dimension $n_1\geq2$ and $X_2$ any compact complex manifold; put $n=n_1+\dim X_2$. Then $X_1\times X_2$ is divisorially hyperbolic.
\end{Lem}

\begin{proof}
Let $f=(f_1,f_2):\C^{n-1}\to X_1\times X_2$ be holomorphic. The restriction of $f_1$ to each complex line is an entire curve in $X_1$, hence constant, and two points of $\C^{n-1}$ lie on a common line: $f_1$ is constant. Then
\begin{eqnarray*}
\mathrm{rank}\,df & \leq & \dim X_2\;=\;n-n_1\;\leq\;n-2\;<\;n-1,
\end{eqnarray*}
so no $f$ is non-degenerate at a point.
\end{proof}

We use three classical facts.
\begin{enumerate}[label=(F\arabic*)]
\item {\bf Lefschetz hyperplane theorem.} Let $M$ be a smooth complex projective variety and $D\subset M$ a smooth ample hypersurface with $\dim D\geq2$. Then $\pi_1(D)\simeq\pi_1(M)$. In particular a smooth hypersurface of dimension $\geq2$ of $\mathbb P^N$ is simply connected.
\item {\bf Bloch--Ochiai--Green} [Och77, Gre75]. The Zariski closure of the image of a non-constant entire curve $\C\to A$ in an abelian variety $A$ is a translate $B+a$ of an abelian subvariety $B\neq0$. An abelian variety is {\it simple} if its only abelian subvarieties are $0$ and $A$ (a very general abelian variety is simple). Consequently every proper closed subvariety of a simple abelian variety is Brody hyperbolic.
\item {\bf Hyperbolic hypersurfaces.} Smooth Brody (equivalently Kobayashi) hyperbolic surfaces of $\mathbb P^3$ exist, e.g.\ of degree $6$ [Duv04] and, very generically, of every degree $\geq21$ [DEG00]; smooth hyperbolic hypersurfaces of $\mathbb P^{N}$ exist for every $N\geq3$ in sufficiently large degree [MN96].
\end{enumerate}

In each example below we check, for a compact K\"ahler (indeed projective) manifold $X$: (a) that $X$ is Brody hyperbolic, hence divisorially hyperbolic (Section~\ref{sec:prelim}); (b) the topological data (fundamental group); (c) the obstruction showing that $X$ is neither balanced nor Gauduchon hyperbolic. When $\pi_1(X)$ is trivial the obstruction is elementary: then $\widetilde X=X$ and every form on $X$ is bounded, so balanced (resp.\ Gauduchon) hyperbolicity would give $\omega^{n-1}=d\Gamma$ (resp.\ $\partial\Gamma+\bar\partial\Gamma'$) on the compact manifold $X$, which contradicts $\int_X\omega^{n-1}\wedge\omega_K>0$ for a K\"ahler form $\omega_K$ (Stokes).

\begin{Ex}[surfaces: $n=2$]\label{ex:surfaces}
For $n=2$, divisorially hyperbolic means Brody hyperbolic (Proposition~\ref{prop:brody}), and balanced (resp.\ Gauduchon) hyperbolic means K\"ahler (resp.\ SKT) hyperbolic (Section~\ref{sec:prelim}). So the examples below are Brody hyperbolic surfaces that are not K\"ahler hyperbolic; they also show that the implication SKT hyperbolic $\Rightarrow$ Kobayashi hyperbolic is strict.
\begin{enumerate}[label=(\arabic*)]
\item {\it Hyperbolic surfaces of $\mathbb P^3$.} Let $S\subset\mathbb P^3$ be a smooth hyperbolic surface (F3).
 \begin{enumerate}[label=(\alph*)]
 \item $S$ is projective, hence K\"ahler, and Brody hyperbolic, so divisorially hyperbolic.
 \item $S$ is simply connected by (F1).
 \item Since $\pi_1(S)=1$ is amenable, $S$ is neither K\"ahler nor SKT hyperbolic (Corollary~\ref{cor:amen}); concretely, a K\"ahler hyperbolic metric $\omega$ would satisfy $\omega=d\eta$ on the compact surface $S$, contradicting $\int_S\omega^2>0$.
 \end{enumerate}
\item {\it Ample divisors of simple abelian threefolds.} Let $A=\C^3/\Lambda$ be a simple abelian threefold, $L$ a very ample line bundle on $A$ and $D\in|L|$ a smooth member (Bertini); $D$ is a smooth ample divisor, a projective surface.
 \begin{enumerate}[label=(\alph*)]
 \item $D$ is Brody hyperbolic: if $f:\C\to D\subset A$ is non-constant, by (F2) the closure of $f(\C)$ is a translate $B+a\subset D$ of an abelian subvariety $B\neq0$; as $A$ is simple $B=A$, which is impossible since $D\neq A$.
 \item $\pi_1(D)\simeq\pi_1(A)=\Z^6$ by (F1).
 \item $\Z^6$ is amenable (equivalently, the universal cover has polynomial volume growth, so its Cheeger constant is $0$), so $D$ is not K\"ahler hyperbolic by Corollary~\ref{cor:amen}. This gives a Brody hyperbolic surface with infinite fundamental group that is not K\"ahler hyperbolic.
 \end{enumerate}
\end{enumerate}
\end{Ex}

\begin{Ex}[higher dimensions: $n\geq3$]\label{ex:higher}
\begin{enumerate}[label=(\arabic*)]
\item {\it Hyperbolic hypersurfaces of $\mathbb P^{n+1}$.} Let $X\subset\mathbb P^{n+1}$ be a smooth hyperbolic hypersurface (F3). It is projective and Brody hyperbolic, hence divisorially hyperbolic; it is simply connected by (F1) ($\dim X=n\geq2$); so it is neither balanced nor Gauduchon hyperbolic by the elementary argument above.
\item {\it Ample divisors of simple abelian varieties.} Let $A$ be a simple abelian variety of dimension $n+1$ and $D\subset A$ a smooth ample divisor (Bertini). Exactly as in Example~\ref{ex:surfaces}(2), $D$ is Brody hyperbolic by (F2) and $\pi_1(D)\simeq\Z^{2n+2}$ by (F1) ($\dim D=n\geq2$). Hence $D$ is divisorially hyperbolic, and, $\Z^{2n+2}$ being amenable, neither balanced nor Gauduchon hyperbolic by Corollary~\ref{cor:amen}.
\item {\it Products with a curve of genus $\geq2$.} Let $X_1$ be as in (1) (of dimension $n_1\geq2$) and $Y$ a compact Riemann surface of genus $g\geq2$; put $X=X_1\times Y$, $n=n_1+1\geq3$.
 \begin{enumerate}[label=(\alph*)]
 \item $X$ is projective and Brody hyperbolic (a product of Brody hyperbolic manifolds), hence divisorially hyperbolic.
 \item $\pi_1(X)=\pi_1(X_1)\times\pi_1(Y)=\pi_1(Y)$ is {\it non-amenable}, so Corollary~\ref{cor:amen} gives no information.
 \item The map $f:X_1\to X$, $x\mapsto(x,y_0)$, is an embedding of a compact simply connected manifold of dimension $n-1$, so $df$ has rank $n-1$ everywhere, and Proposition~\ref{prop:hyper} shows that $X$ is neither balanced nor Gauduchon hyperbolic. Geometrically, the universal cover is $X_1\times\mathbb D$ and contains the compact hypersurface $X_1\times\{pt\}$, on which $\omega^{n-1}$ would integrate to $0$.
 \end{enumerate}
 This also shows that non-amenability of $\pi_1$ is not sufficient for balanced hyperbolicity in dimension $\geq3$.
\item {\it Products with a projective space.} Let $X_1$ be as in (1) and $X=X_1\times\mathbb P^m$, $m\geq1$. By Lemma~\ref{lem:prod} $X$ is divisorially hyperbolic. It is simply connected and projective, hence neither balanced nor Gauduchon hyperbolic by the elementary argument above. It is not Brody hyperbolic ($\mathbb P^m\subset X$): in dimension $\geq3$ divisorial hyperbolicity does not imply Brody hyperbolicity, in contrast with Proposition~\ref{prop:brody}.
\end{enumerate}
\end{Ex}

\begin{Question}
Is a Brody hyperbolic compact K\"ahler surface with non-amenable fundamental group K\"ahler hyperbolic? (The answer is no in dimension $\geq3$: Example~\ref{ex:higher}(3).)
\end{Question}

\begin{table}[h]
\centering
\footnotesize
\setlength{\tabcolsep}{3pt}
\begin{tabular}{|p{0.31\linewidth}|p{0.07\linewidth}|p{0.14\linewidth}|p{0.09\linewidth}|p{0.27\linewidth}|}
\hline
Example & $n$ & $\pi_1(X)$ & Brody hyp. & Obstruction to balanced and Gauduchon hyp.\\
\hline
Hyperbolic surface of $\mathbb P^3$ & $2$ & $1$ & yes & amenable $\pi_1$\\
\hline
Ample divisor of a simple abelian $3$-fold & $2$ & $\Z^6$ & yes & amenable $\pi_1$\\
\hline
Hyperbolic hypersurface of $\mathbb P^{n+1}$ & $\geq3$ & $1$ & yes & amenable $\pi_1$\\
\hline
Ample divisor of a simple abelian $(n+1)$-fold & $\geq3$ & $\Z^{2n+2}$ & yes & amenable $\pi_1$\\
\hline
$X_1\times Y$, $g(Y)\geq2$ & $\geq3$ & non-amenable & yes & liftable hypersurface $X_1\times\{y_0\}$\\
\hline
$X_1\times\mathbb P^m$ & $\geq3$ & $1$ & no & amenable $\pi_1$\\
\hline
\end{tabular}
\caption{Divisorially hyperbolic manifolds that are neither balanced nor Gauduchon hyperbolic.}\label{tab:ex}
\end{table}


\begin{thebibliography}{99}
\bibitem[BHPV04]{BHPV04} W.~Barth, K.~Hulek, C.~Peters, A.~Van de Ven, {\it Compact complex surfaces}, 2nd ed., Springer, 2004.
\bibitem[Bro78]{Bro78} R.~Brody, {\it Compact manifolds and hyperbolicity}, Trans.\ Amer.\ Math.\ Soc.\ {\bf 235} (1978), 213--219.
\bibitem[Bro81]{Bro81} R.~Brooks, {\it The fundamental group and the spectrum of the Laplacian}, Comment.\ Math.\ Helv.\ {\bf 56} (1981), 581--598.
\bibitem[DEG00]{DEG00} J.-P.~Demailly, J.~El Goul, {\it Hyperbolicity of generic surfaces of high degree in projective 3-space}, Amer.\ J.\ Math.\ {\bf 122} (2000), 515--546.
\bibitem[DP04]{DP04} J.-P.~Demailly, M.~P\u{a}un, {\it Numerical characterization of the K\"ahler cone of a compact K\"ahler manifold}, Ann.\ of Math.\ {\bf 159} (2004), 1247--1274.
\bibitem[Duv04]{Duv04} J.~Duval, {\it Une sextique hyperbolique dans $\mathbb P^3(\mathbb C)$}, Math.\ Ann.\ {\bf 330} (2004), 473--476.
\bibitem[Gre75]{Gre75} M.~Green, {\it Some Picard theorems for holomorphic maps to algebraic varieties}, Amer.\ J.\ Math.\ {\bf 97} (1975), 43--75.
\bibitem[Gro91]{Gro91} M.~Gromov, {\it K\"ahler hyperbolicity and $L_2$-Hodge theory}, J.\ Differential Geom.\ {\bf 33} (1991), 263--292.
\bibitem[Khe26]{Khe26} A.~Khelifati, {\it Deformations of non-K\"ahler hyperbolicity notions and modifications of degenerate balanced manifolds}, arXiv:2512.19284.
\bibitem[Mar24]{Mar24} S.~Marouani, {\it SKT hyperbolic and Gauduchon hyperbolic compact complex manifolds}, Bull.\ Soc.\ Math.\ France {\bf 152} (2024), no.~3, 519--545; arXiv:2305.08122.
\bibitem[MN96]{MN96} K.~Masuda, J.~Noguchi, {\it A construction of hyperbolic hypersurface of $\mathbb P^n(\mathbb C)$}, Math.\ Ann.\ {\bf 304} (1996), 339--362.
\bibitem[MP22]{MP22} S.~Marouani, D.~Popovici, {\it Balanced hyperbolic and divisorially hyperbolic compact complex manifolds}, Math.\ Res.\ Lett.\ {\bf 30} (2023), no.~6, 1813--1855.
\bibitem[Och77]{Och77} T.~Ochiai, {\it On holomorphic curves in algebraic varieties with ample irregularity}, Invent.\ Math.\ {\bf 43} (1977), 83--96.
\bibitem[ST10]{ST10} J.~Streets, G.~Tian, {\it A parabolic flow of pluriclosed metrics}, Int.\ Math.\ Res.\ Not.\ IMRN 2010, no.~16, 3101--3133.
\end{thebibliography}
\end{document}